\documentclass[preprint,twosides,times]{elsarticle}
\usepackage[utf8]{inputenc}
\usepackage{bussproofs}
\usepackage{array}
\usepackage{amsthm}
\usepackage{amssymb}
\usepackage{amsmath}
\usepackage[all]{xy}
\usepackage{stmaryrd}
\usepackage{float}
\usepackage{graphicx}
\floatstyle{boxed}
\restylefloat{figure}
\restylefloat{table}
\usepackage{hyperref}
\hypersetup{colorlinks=false}

\swapnumbers
\theoremstyle{definition}
\newtheorem{definition}{Definition}[section]
\newtheorem*{definition*}{Definition}

\theoremstyle{plain}
\newtheorem{proposition}[definition]{Proposition}
\newtheorem{lemma}[definition]{Lemma}
\newtheorem{theorem}[definition]{Theorem}
\newtheorem{corollary}[definition]{Corollary}
\newtheorem*{lemma*}{Lemma}
\newtheorem*{theorem*}{Theorem}
\newtheorem*{corollary*}{Corollary}

\theoremstyle{remark}
\newtheorem{remark}[definition]{Remark}
\newtheorem{example}[definition]{Example}

\newcommand{\axc}[1]{\AxiomC{#1}}
\newcommand{\uic}[2]{\RightLabel{\small{#2}}\UnaryInfC{#1}}
\newcommand{\bic}[2]{\RightLabel{\small{#2}}\BinaryInfC{#1}}
\newcommand{\tic}[2]{\RightLabel{\small{#2}}\TrinaryInfC{#1}}

\newcommand{\mqcplus}{MQC$^+$}
\newcommand{\dnst}{DNS$_T$}
\newcommand{\dnstforall}{DNS$^\forall_T$}
\newcommand{\dnstimpl}{DNS$^\Rightarrow_T$}
\newcommand{\mpt}{MP$_T$}
\newcommand{\cchoice}{$\text{AC}_0$}
\newcommand{\cchoicen}{$\text{AC}_0^\text{N}$}
\newcommand{\dchoice}{$\text{DC}$}

\newcommand{\caseof}[5]{\mathsf{case} ~ #1 ~\mathsf{of} ~ \left(#2.#3 \| #4.#5\right)}
\newcommand{\destas}[4]{\mathsf{dest} ~ #1 ~\mathsf{as} ~ \left(#2.#3\right) ~\mathsf{in}~ #4}
\newcommand{\shift}[2]{\mathcal{S}#1.#2}
\newcommand{\reset}[1]{\# #1}

\begin{document}

\begin{frontmatter}
  \title{Delimited control operators prove Double-negation Shift}
  \author{Danko Ilik\footnote{Present address: Faculty of Informatics, University ``Goce Delčev'', PO Box 201, 2000 Štip, Macedonia; E-mail: danko.ilik@ugd.edu.mk}}

  \address{Ecole Polytechnique, INRIA, CNRS \& Université Paris Diderot\\
    Address: INRIA PI-R2, 23 avenue d'Italie, CS 81321, 75214 Paris Cedex 13, France\\
    E-mail: danko.ilik@polytechnique.edu}

  \begin{abstract}
    We propose an extension of minimal intuitionistic predicate logic, based on delimited control operators, that can derive the predicate-logic version of the Double-negation Shift schema, while preserving the disjunction and existence properties.
  \end{abstract}


\end{frontmatter}

\section{Introduction}
\label{delcont}
In \cite{Herbelin2010}, Hugo Herbelin showed that, by extending the proof-term calculus of intuitionistic predicate logic with a restricted form of \emph{delimited control operators}, one can obtain a logical system able to derive a predicate-logic version of Markov's Principle, $\neg\neg\exists x A(x) \Rightarrow \exists x A(x)$ (for $A(x)$ a $\{\Rightarrow,\forall\}$-free formula), while remaining essentially intuitionistic -- satisfying the disjunction and existence properties. 

Separately, \cite{HerbelinPC} he also observed that using the full power of delimited control operators one can derive the predicate-logic version of the Double-negation Shift schema, $\forall x \neg\neg A(x) \Rightarrow \neg\neg\forall x A(x)$ (where $A(x)$ is arbitrary), and posed the question whether there is a corresponding logical system which also posseses the disjunction and existence properties. With this article, we answer Herbelin's question in the affirmative.

Delimited control operators have appeared in Theoretical Computer Science, in Semantics of Programming Languages, as a powerful abstraction to account for so-called \emph{computational effects}. While being pervasive in the practise of writing computer programs (for, they include facilities as basic as reading from and writing into memory, stopping the execution of the program, or parallel computation), giving a good mathematical explanation of effects is still one of the major research topics in Semantics. 

An important step in that direction was a result of Filinski \cite{Filinski1994,FilinskiThesis}, who showed that every monadic computational effect can be operationally simulated by the delimited control operators \emph{shift/reset}, introduced previously by himself and Danvy \cite{DanvyF1989,DanvyF1990}. However, the logical status of shift/reset themselves, and other such operators in general, remains to be fully established, something we hope to being contributing to with this article.

Another interest in delimited control operators, and actually our original interest in them, comes from the role they promise to be playing in a future constructive proof of completeness of full intuitionistic logic (with $\vee$ and $\exists$) with respect to Kripke semantics, something described in \cite{IlikCPS,Ilik}.

We illustrate the utility of delimited control operators, by considering some examples in a $\lambda$-calculus extended with them, containing also natural numbers with the plus operation. The extension consists of two constructs, a \emph{delimiter} ($\#$ -- ``reset'')  and a \emph{control operator} ($\mathcal{S}$ -- ``shift''). The delimiter is used as a special kind of brackets in a $\lambda$-term, so that the control operator, which can only appear inside such ``brackets'', be able to gain control of its surrounding context, up to the delimiter. For example, in the following $\lambda$-term reduction,
\[
\begin{array}{lllllll}
  1 + \reset{2 + \shift{k}{4}} &\to& 1 + \reset{4\left\{\left.(\lambda a.\reset{2+a})\right/k\right\}} &=& 1+\reset{4} &\to& 1+4 \to 5,
\end{array}
\]
reset is used to delimit the sub-term $2 + \shift{k}{4}$. Shift then behaves as a binder, alike $\lambda$-abstraction, that names the abstracted surroundings of shift, $2 + \Box$, by $k$, and replaces in its sub-expression, $4$, all occurrences of $k$ by the abstracted surroundings. In this case, $k$ is not used inside shift -- this corresponds to the so-called ``exceptions'' effect that Herbelin found out to be the computational contents behind Markov's Principle. In the next example, $k$ is used; the sub-term inside shift uses its surrounding context twice:
\begin{align*}
        & 1+\reset{2+\shift{k}{k 4 + k 8}} \\
  \to   & 1+\reset{(\lambda a.\reset{2+a}) 4 + (\lambda a.\reset{2+a}) 8} \\
  \to^+ & 1+\reset{(\reset{6}) + (\reset{10})}\\
  \to^+ & 1+\reset{6 + 10}\\
  \to^+ & 17
\end{align*}


From the logical perspective, considering natural deduction formalisms which can be isomorphically presented by proof-$\lambda$-terms, we see delimited control, when added to the syntax of such proof terms, as a means of being able to access \emph{a certain part} of the surroundings of a proof term from inside the proof term itself.\footnote{This is to be contrasted to what happens with the (undelimited) control operator \textit{call/cc}, which is better known in Logic for its role in the development of classical realisability \cite{Griffin1990, Krivine2001, Krivine2010a, Krivine2010b} -- call/cc amounts computationally to aborting the entire computation and, since its effect is not delimited, one has no hope of getting a natural computational interpretation from classical realisability: a realiser of an existential statement needs not be a program which computes a witness for the existential quantifier.} The part of the surrounding that we want to be able to access will be defined as a ``pure evaluation context'' in Section~\ref{delcont_mqcplus}; logically, it is the surroundings of a proof term for a $\{\Rightarrow,\forall\}$-free formula,\footnote{Following Berger \cite{Berger2004}, we call the $\{\Rightarrow,\forall\}$-free formulae, $\Sigma$-formulae, and denote them by $S, T, U$, while general formulae are denoted by $A, B, C$.} which is the predicate logic equivalent of arithmetic $\Sigma^0_1$-formulae, for which we know that classical and intuitionistic provability coincide. In other words, we propose a proof-term calculus for a logic which is essentially intuitionistic, except that at the fragment ``$\Sigma^0_1$'' we are allowed to use classical reasoning to obtain more (succinct) proofs.

The paper is organised as follows. In the next Section~\ref{delcont_mqcplus}, we introduce our system \mqcplus. The acronym comes from Troelstra: IQC is intuitionistic predicate logic, MQC is minimal predicate logic (IQC without the $\bot_E$ rule), and CQC is classical predicate logic. In Section~\ref{delcont_equicons}, we characterise the relationship between \mqcplus, MQC, and CQC; in particular, we show that an extension to \emph{predicate} logic of Glivenko's Theorem holds for our system, unlike for MQC. In Section~\ref{delcont_srprogress}, we prove properties of the reduction relation on proof terms, from which we obtain the  Disjunction and Existence Property for closed derivations of \mqcplus. In the final Section~\ref{delcont_future}, we discuss related and future work.


\section{The system \mqcplus}
\label{delcont_mqcplus}
The natural deduction system of \mqcplus\ is shown in Table~\ref{tab:mqcplus}. It consists of the proof rules of minimal intuitionistic predicate logic MQC, plus two new ones, ``shift'' $(\mathcal{S})$ and ``reset'' $(\#)$. 

\begin{table}
  \centering
  \begin{tabular}{ m{3cm} m{8cm} }
    \multicolumn{2}{ m{11cm} }{
      \begin{prooftree}
        \axc{$A \in \Gamma$}
        \uic{$\Gamma\vdash_\diamond A$}{\textsc{Ax}}
      \end{prooftree}
    }
    \\
    \begin{prooftree}
      \axc{$\Gamma\vdash_\diamond A_1$}
      \axc{$\Gamma\vdash_\diamond A_2$}
      \bic{$\Gamma\vdash_\diamond A_1\wedge A_2$}{$\wedge_I$}
    \end{prooftree}
    &
    \begin{prooftree}
      \axc{$\Gamma\vdash_\diamond A_1\wedge A_2$}
      \uic{$\Gamma\vdash_\diamond A_i$}{$\wedge^i_E$}
    \end{prooftree}
    \\
    \begin{prooftree}
      \axc{$\Gamma\vdash_\diamond A_i$}
      \uic{$\Gamma\vdash_\diamond A_1\vee A_2$}{$\vee^i_I$}
    \end{prooftree}
    &
    \begin{prooftree}
      \axc{$\Gamma\vdash_\diamond A_1\vee A_2$}
      \axc{$\Gamma, A_1\vdash_\diamond C$}
      \axc{$\Gamma, A_2\vdash_\diamond C$}
      \tic{$\Gamma\vdash_\diamond C$}{$\vee_E$}
    \end{prooftree}
    \\
    \begin{prooftree}
      \axc{$\Gamma, A_1\vdash_\diamond A_2$}
      \uic{$\Gamma\vdash_\diamond A_1\Rightarrow A_2$}{$\Rightarrow_I$}
    \end{prooftree}
    &
    \begin{prooftree}
      \axc{$\Gamma\vdash_\diamond A_1\Rightarrow A_2$}
      \axc{$\Gamma\vdash_\diamond A_1$}
      \bic{$\Gamma\vdash_\diamond A_2$}{$\Rightarrow_E$}
    \end{prooftree}
    \\
    \begin{prooftree}
      \axc{$\Gamma\vdash_\diamond A(x)$}
      \axc{$x\text{-fresh}$}
      \bic{$\Gamma\vdash_\diamond \forall x A(x)$}{$\forall_I$}
    \end{prooftree}
    &
    \begin{prooftree}
      \axc{$\Gamma\vdash_\diamond \forall x A(x)$}
      \uic{$\Gamma\vdash_\diamond A(t)$}{$\forall_E$}
    \end{prooftree}
    \\
    \begin{prooftree}
      \axc{$\Gamma\vdash_\diamond A(t)$}
      \uic{$\Gamma\vdash_\diamond \exists x.A(x)$}{$\exists_I$}
    \end{prooftree}
    &
    \begin{prooftree}
      \axc{$\Gamma\vdash_\diamond \exists x.A(x)$}
      \axc{$\Gamma, A(x)\vdash_\diamond C$}
      \axc{$x\text{-fresh}$}
      \tic{$\Gamma\vdash_\diamond C$}{$\exists_E$}
    \end{prooftree}
    \\
    \begin{prooftree}
      \axc{$\Gamma\vdash_T T$}
      \uic{$\Gamma\vdash_\diamond T$}{$\#$ (``reset'')}
    \end{prooftree}
    &
    \begin{prooftree}
      \axc{$\Gamma, A\Rightarrow T\vdash_T T$}
      \uic{$\Gamma\vdash_T A$}{$\mathcal{S}$ (``shift'')}
    \end{prooftree}\\
  \end{tabular}
  
  \caption[\mqcplus]{Natural deduction system of \mqcplus}
  \label{tab:mqcplus}
\end{table}

The turnstile symbol ``$\vdash$'' can carry an annotation -- a $\Sigma$-formula $T$ -- which is neither used nor changed by the intuitionistic rules. We use the wild-card symbol $\diamond$ for this purpose, to mean that there either is an annotating formula $T$, or that there is none. In the proof rules where the wild-card appears both above and below the line, it means that either there is the same annotation both above and below, or that there is no annotation above and no annotation below. 

The rule $(\#)$ can only be applied when the conclusion is a $\Sigma$-formula $T$. It acts as a delimiter in the proof tree, (re-)initialising the annotation with the formula $T$; from that point upwards in the tree, classical reasoning is allowed -- but, only so because we are ultimately proving a $\Sigma$-formula. The rule $(\mathcal{S})$ can then be used, inside a sub-tree with $(\#)$ at its root, as a kind of $(\neg\neg_E)$ rule. Its role is to ``escape'' to the nearest enclosing delimiter once a witness for the $\Sigma$-formula from the annotation has been found.

However, note that, although there can be arbitrarily many uses of the $(\#)$ and $(\mathcal{S})$ rules in a derivation tree, \emph{only one} formula $T$ is allowed to appear in annotations, globally, of a derivation tree. This means that the global $T$ is set once and for all, hence it is not possible to mix derivations using different $T$ and $T'$. Were we in IQC, a natural choice for the global $T$ would have been $\bot$.

As examples, we give the derivations for (generalisations of) the minimal-predicate-logic versions\footnote{The distinguished formula $T$ plays the role of $\bot$ and the hypothesis $T\Rightarrow S$ plays the role of the $\bot_E$ rule.} of Markov's Principle,
\begin{equation}\tag{\mpt}
(T\Rightarrow S) \Rightarrow ((S\Rightarrow T)\Rightarrow T) \Rightarrow S,
\end{equation}\label{example_derivations}
and Double-negation Shift,
\begin{equation}\tag{\dnst}
  \forall x \left(\left(A(x)\Rightarrow T\right)\Rightarrow T\right)
  \Rightarrow \left(\forall x A(x) \Rightarrow T\right) \Rightarrow T,
\end{equation}
where, according to the already set convention, $T$ and $S$ are $\Sigma$-formulae, while $A(x)$ is a general one.

\begin{prooftree}
\axc{}
\uic{$\cdots \vdash_S T\Rightarrow S$}{\textsc{Ax}}
\axc{}
\uic{$\cdots \vdash_S (S\Rightarrow T)\Rightarrow T$}{\textsc{Ax}}
\axc{}
\uic{$\cdots, S \vdash_S S$}{$\textsc{Ax}$}
\uic{$\cdots, S \vdash_S T$}{$\mathcal{S}$}
\uic{$\cdots \vdash_S S\Rightarrow T$}{$\Rightarrow_I$}
\bic{$\cdots \vdash_S T$}{$\Rightarrow_E$}
\bic{$T\Rightarrow S, (S\Rightarrow T)\Rightarrow T \vdash_S S$}{$\Rightarrow_E$}
\uic{$T\Rightarrow S, (S\Rightarrow T)\Rightarrow T \vdash S$}{$\#$}
\uic{$T\Rightarrow S \vdash ((S\Rightarrow T)\Rightarrow T) \Rightarrow S$}{$\Rightarrow_I$}
\uic{$\vdash (T\Rightarrow S) \Rightarrow ((S\Rightarrow T)\Rightarrow T) \Rightarrow S$}{$\Rightarrow_I$}
\end{prooftree}

~

\begin{prooftree}
\axc{}
\uic{$\cdots\vdash\forall x A(x)\Rightarrow T$}{\textsc{Ax}}
\axc{$\vdots$}
\uic{$\cdots, \forall x  \left(\left(A(x)\Rightarrow T\right)\Rightarrow T\right), A(x)\Rightarrow T \vdash_T T$}{$\forall_E, \Rightarrow_E$, \textsc{Ax}}
\uic{$\cdots, \forall x  \left(\left(A(x)\Rightarrow T\right)\Rightarrow T\right) \vdash_T A(x)$}{$\mathcal{S}$}
\uic{$\cdots, \forall x  \left(\left(A(x)\Rightarrow T\right)\Rightarrow T\right) \vdash_T \forall x  A(x)$}{$\forall_I$, $x$-fresh}
\bic{$\forall x A(x) \Rightarrow T, \forall x  \left(\left(A(x)\Rightarrow T\right)\Rightarrow T\right)\vdash_T T$}{$\Rightarrow_E$}
\uic{$\forall x A(x) \Rightarrow T, \forall x \left(\left(A(x)\Rightarrow T\right)\Rightarrow T\right) \vdash T$}{$\#$}
\uic{$\forall x  \left(\left(A(x)\Rightarrow T\right)\Rightarrow T\right) \vdash \left(\forall x  A(x) \Rightarrow T\right) \Rightarrow T$}{$\Rightarrow_I$}
\uic{$\vdash \forall x  \left(\left(A(x)\Rightarrow T\right)\Rightarrow T\right) \Rightarrow \left(\forall x  A(x) \Rightarrow T\right) \Rightarrow T$}{$\Rightarrow_I$}
\end{prooftree}

~

We now define a calculus of proof-term annotations for the natural deduction system of \mqcplus, a version of simply typed $\lambda$-calculus with constants for handling all logical connectives and the delimited control operators, and then a reduction system for proof terms; the idea is that reducing a proof term describes the process of normalising a natural deduction derivation.

 The definitions are based on standard treatments of Logic as $\lambda$-calculus (see, for example, \cite{SorensenU}), and standard treatment of $\lambda$-calculus with shift/reset from Semantics of Programming Languages (for example, \cite{AsaiK2007}). 

\begin{table}
  \centering
  \begin{tabular}{ m{4cm} m{7cm} }
    \multicolumn{2}{ m{11cm} }{
      \begin{prooftree}
        \axc{$(a:A) \in \Gamma$}
        \uic{$\Gamma\vdash_\diamond a:A$}{\textsc{Ax}}
      \end{prooftree}
    }
    \\
    \begin{prooftree}
      \axc{$\Gamma\vdash_\diamond p:A_1$}
      \axc{$\Gamma\vdash_\diamond q:A_2$}
      \bic{$\Gamma\vdash_\diamond (p,q):A_1\wedge A_2$}{$\wedge_I$}
    \end{prooftree}
    &
    \begin{prooftree}
      \axc{$\Gamma\vdash_\diamond p:A_1\wedge A_2$}
      \uic{$\Gamma\vdash_\diamond \pi_i p:A_i$}{$\wedge^i_E$}
    \end{prooftree}
    \\
    \begin{prooftree}
      \axc{$\Gamma\vdash_\diamond p:A_i$}
      \uic{$\Gamma\vdash_\diamond \iota_i p:A_1\vee A_2$}{$\vee^i_I$}
    \end{prooftree}
    & \\
    \multicolumn{2}{ m{11cm} }{
      \begin{prooftree}
        \axc{$\Gamma\vdash_\diamond p:A_1\vee A_2$}
        \axc{$\Gamma, a_1:A_1\vdash_\diamond q_1:C$}
        \axc{$\Gamma, a_2:A_2\vdash_\diamond q_2:C$}
        \tic{$\Gamma\vdash_\diamond \caseof{p}{a_1}{q_1}{a_2}{q_2}:C$}{$\vee_E$}
      \end{prooftree}
    }
    \\
    \begin{prooftree}
      \axc{$\Gamma, a:A_1\vdash_\diamond p:A_2$}
      \uic{$\Gamma\vdash_\diamond \lambda a.p:A_1\Rightarrow A_2$}{$\Rightarrow_I$}
    \end{prooftree}
    &
    \begin{prooftree}
      \axc{$\Gamma\vdash_\diamond p:A_1\Rightarrow A_2$}
      \axc{$\Gamma\vdash_\diamond q:A_1$}
      \bic{$\Gamma\vdash_\diamond p q:A_2$}{$\Rightarrow_E$}
    \end{prooftree}
    \\
    \begin{prooftree}
      \axc{$\Gamma\vdash_\diamond p:A(x)$}
      \axc{$x\text{-fresh}$}
      \bic{$\Gamma\vdash_\diamond \lambda x.p:\forall x  A(x)$}{$\forall_I$}
    \end{prooftree}
    &
    \begin{prooftree}
      \axc{$\Gamma\vdash_\diamond p:\forall x A(x)$}
      \uic{$\Gamma\vdash_\diamond p t:A(t)$}{$\forall_E$}
    \end{prooftree}
    \\
    \begin{prooftree}
      \axc{$\Gamma\vdash_\diamond p:A(t)$}
      \uic{$\Gamma\vdash_\diamond (t,p):\exists x.A(x)$}{$\exists_I$}
    \end{prooftree}
    & \\
    \multicolumn{2}{ m{11cm} }{
      \begin{prooftree}
        \axc{$\Gamma\vdash_\diamond p:\exists x.A(x)$}
        \axc{$\Gamma, a:A(x)\vdash_\diamond q:C$}
        \axc{$x\text{-fresh}$}
        \tic{$\Gamma\vdash_\diamond \destas{p}{x}{a}{q}:C$}{$\exists_E$}
      \end{prooftree}
    }
    \\
    \begin{prooftree}
      \axc{$\Gamma\vdash_T p:T$}
      \uic{$\Gamma\vdash_\diamond \reset{p}:T$}{$\#$ (``reset'')}
    \end{prooftree}
    &
    \begin{prooftree}
      \axc{$\Gamma, k:A\Rightarrow T\vdash_T p:T$}
      \uic{$\Gamma\vdash_T \shift{k}{p}:A$}{$\mathcal{S}$ (``shift'')}
    \end{prooftree}\\
  \end{tabular}
  
  \caption[\mqcplus\ with proof terms]{Proof term annotation for the natural deduction system of \mqcplus}
  \label{tab:mqcplus_terms}
\end{table}

\begin{definition} The set of \emph{proof terms} is defined by the following inductive definition,
  \begin{multline*}
    p,q ::= a ~|~ \iota_1 p ~|~ \iota_2 p ~|~ \caseof{p}{a_1}{q_1}{a_2}{q_2} ~|~ (p,q) ~|~ \pi_1 p ~|~ \pi_2 p ~|~ \lambda a.p ~|~ p q ~|~\\
    \lambda x.p ~|~ p t ~|~ (t,p) ~|~ \destas{p}{x}{a}{q} ~|~ \reset{p} ~|~ \shift{k}{p}
  \end{multline*}
  where $a,b,k,l$ denote hypothesis variables, $x,y,z$ denote quantifier variables, and $t,u,v$ denote quantifier terms (individuals); hence, $\lambda a.p$ is a constructor for implication, while $\lambda x.p$ is a constructor for universal quantification; $(p,q)$ is a constructor for conjunction while $(t,p)$ is a constructor for existential quantification, and $p q$ is a destructor for implication while $p t$ is a destructor for universal quantification.

\end{definition}

\begin{remark}The $\mathcal{S}$ in $\shift{k}{p}$ is a binder, it binds $k$ in $p$ just as $\lambda$ binds $a$ in $q$ in a lambda abstraction $\lambda a.q$. Following standard terminology, we sometimes call $k$ a \textit{continuation} variable.
\end{remark}

\begin{definition}\label{values} The subset of proof terms known as \emph{values} is defined by:
  \[
    V ::= a ~|~ \iota_1 V ~|~ \iota_2 V ~|~ (V,V) ~|~ (t,V) ~|~ \lambda a.p ~|~ \lambda x.p
  \]
\end{definition}

\begin{definition}\label{pureecontext} The set of \emph{pure evaluation contexts}, a subset of all proof terms with one placeholder or ``hole'', is defined by:
  \begin{multline*}
    P ::= [~] ~|~ \caseof{P}{a_1}{p_1}{a_2}{p_2} ~|~ \pi_1P ~|~
    \pi_2P ~|~ \destas{P}{x}{a}{p} ~|~ \\
    Pq ~|~ (\lambda a.q)P ~|~ Pt ~|~ \iota_1P ~|~ \iota_2P ~|~ (P,p)
    ~|~ (V,P) ~|~ (t,P)
  \end{multline*}
\end{definition}

The association of proof terms to natural deduction derivations is given in Table~\ref{tab:mqcplus_terms}. $P[p]$ denotes the proof term obtained from $P$ by replacing its placeholder $[~]$ with the proof term $p$.

In order to define a reduction relation on proof terms we also need the notion of (non-pure) evaluation context.

\begin{definition} The set of \emph{evaluation contexts} is given by the following inductive definition:
  \begin{multline*}
    E ::= [~] ~|~ \caseof{E}{a_1}{p_1}{a_2}{p_2} ~|~ \pi_1E ~|~
    \pi_2E ~|~ \destas{E}{x}{a}{p} ~|~ \\
    Eq ~|~ (\lambda a.q)E ~|~ Et ~|~ \iota_1E ~|~ \iota_2E
    ~|~ (E,p) ~|~ (V,E) ~|~ (t,E) ~|~ \reset{E}
  \end{multline*}
\end{definition}

The set of evaluation contexts is larger than the set of pure evaluation contexts, because it includes $\#$. As before, $E[p]$ denotes the proof term obtained from $E$ by replacing its placeholder $[~]$ with the proof term $p$.

\begin{definition}\label{mqcplus_reduction} The reduction relation on proof terms ``$\to$'' is defined by the following rewrite rules:
  \begin{align*}
    (\lambda a.p) V &\to p\{V/a\} & \caseof{\iota_i V}{a_1}{p_1}{a_2}{p_2} &\to p_i\{V/a_i\}\\
    (\lambda x.p) t &\to p\{t/x\} & \destas{(t,V)}{x}{a}{p} &\to p\{t/x\}\{V/a\} \\
    \pi_i(V_1,V_2) &\to V_i & \reset{P[\shift{k}{p}]} &\to \reset{p\left\{\left(\lambda a.\reset{P[a]}\right) / k\right\}} \\
    \reset{V} &\to V & E[p]&\to E[p'] \text{ when } p\to p'
  \end{align*}
  The last rule is known as the ``congruent closure'' of the preceding rules. The rule for $\mathcal{S}$ applies only when the evaluation context $P$ is pure. The reduction strategy determined by the rules is standard call-by-value reduction. \cite{Plotkin1975}
\end{definition}

\begin{example}\label{example_derivations_terms} The following are the proof terms corresponding to the derivation trees for \mpt\ and \dnst\ from page \pageref{example_derivations}.
\[
\lambda e. \lambda a. \reset{e (a (\lambda b. \shift{k}{b}))}
\]
\[
\lambda a. \lambda b. \reset{b (\lambda x. \shift{k}{a x k})}
\]
Remark that the proof term for \mpt\ does not make use of the continuation variable $k$, but only uses the $\mathcal{S}$ operator to pass the value $b$, once it has been found in the course of the computation, back to the control delimiter $\#$. 
\end{example}


\section{Relationship to MQC and CQC}
\label{delcont_equicons}
To connect provability in \mqcplus with provability in MQC and CQC, we use the following double-negation translation.


\begin{definition} The \emph{superscript} translation $A^T$ of a
  formula $A$ with respect to a $\Sigma$-formula $T$ is defined via the \emph{subscript} translation $A_T$, which is in turn defined by recursion on the structure of $A$:

  \begin{align*}
    A^T :=& (A_T \Rightarrow T) \Rightarrow T & ~ \\
    ~ & ~ & ~ \\
    A_T :=& A & \text{ if } A \text{ is atomic } \\
    (A\Box B)_T :=& A_T \Box B_T & \text{ for } \Box = \vee,\wedge \\
    (A\Rightarrow B)_T :=& A_T \Rightarrow B^T &  \\
    (\exists A)_T :=& \exists A_T & \\
    (\forall A)_T :=& \forall A^T & \\
  \end{align*}
  We write $\Gamma_T$ for the translation $(-)_T$ applied to each formula of the context $\Gamma$ individually.
\end{definition}

This translation is the standard call-by-value CPS translation of types \cite{Plotkin1975}, and is similar to the Kuroda translation \cite{TroelstraVD1}, the difference being that we add a double negation, not only after $\forall$, but also after $\Rightarrow$. Interestingly, when interpreting, using DNS, the negative translation of the Axiom of Countable Choice \cchoice, a transformation from the Kuroda translation of \cchoice\  into our form, with $\neg\neg$ after $\Rightarrow$, appears to be needed \cite[p. 200]{Kohlenbach}. Also, Avigad has remarked in \cite{AvigadDN} that the Kuroda translation makes essential use of the $\bot_E$ rule.

\begin{remark} When $A$ is a $\Sigma$-formula, we have that $A_T=A$.
\end{remark}

We will denote derivability in \mqcplus\ by ``$\vdash^+$'', derivability in MQC by ``$\vdash^m$'', and the one in CQC by ``$\vdash^c$''. When we say CQC, we have in mind a standard natural deduction calculus, but where $\bot$ is replaced by a distinguished formula $T$ -- which one, will be clear from context -- and correspondingly, the $\bot_E$ rule says that $T\Rightarrow A$, and the $\neg\neg_E$ rule is $(A\Rightarrow T)\Rightarrow T \vdash^c A$.
 The following theorem is not surprising, since, after all, our system is a subsystem of classical logic, but we give it for the sake of completeness, since this version of Kuroda's translation does not use the $\bot_E$ rule in the target system.

\begin{theorem}[Equiconsistency with MQC]\label{equiconsistency}
Given a derivation of $\Gamma\vdash^+ A$, which uses $\mathcal{S}$ and $\mathcal{\#}$ for the $\Sigma$-formula $T$, we can build a derivation of $\Gamma_T\vdash^m A^T$.
\end{theorem}
\begin{proof} By induction on the derivation, using the proof terms listed below. A line above a sub-term marks the place where the induction hypothesis is applied. 
  \begin{align*}
    \overline{a} &= \lambda k. k a \\
    \overline{\lambda a. p} &= \lambda k. k\left(\lambda a.\lambda k'. \overline{p}\left(\lambda b. k' b\right)\right)\\
    \overline{p q} &= \lambda k. \overline{p}\left(\lambda f.\overline{q}\left(\lambda a. f a\left(\lambda b. k b\right)\right)\right)\\
    \overline{(p,q)} &= \lambda k.\overline{p}\left(\lambda a.\overline{q}\left(\lambda b. k\left(a, b\right)\right)\right)\\
    \overline{\pi_1 p} &= \lambda k.\overline{p}\left(\lambda c. k\left(\pi_1 c\right)\right) \\
    \overline{\iota_1 p} &= \lambda k.\overline{p}\left(\lambda a. k\left(\iota_1 a\right)\right)\\
    \overline{\caseof{p}{a_1}{q_1}{a_2}{q_2}} &= \lambda k.\overline{p}\left(\lambda c.~\caseof{ c}{a_1}{\overline{q_1} k}{a_2}{\overline{q_2} k}\right) \\
    \overline{\lambda x. p} &= \lambda k. k\left(\lambda x.\lambda k'. \overline{p}\left(\lambda b. k' b\right)\right)\\
    \overline{p t} &= \lambda k.\overline{p}\left(\lambda f. f t  k\right) \\
    \overline{(t,p)} &= \lambda k.\overline{p}\left(\lambda a. k(t,a)\right) \\
    \overline{\destas{p}{x}{a}{q}} &= \lambda k.\overline{p}\left(\lambda c.~\destas{ c}{x}{a}{\overline{q} k}\right) \\
    \overline{\reset{ a}{p}} &= \lambda k. k\left(\overline{p}(\lambda a.a)\right) \\
    \overline{\shift{l}{p}} &= \lambda k.\left(\overline{p}(\lambda a.a)\right)\left\{\left.\lambda a.\lambda k'.k'\left(k a\right) \right/ l\right\}
  \end{align*}
\end{proof}

In order to relate \mqcplus-provability of certain forms of formulae to their provability in MQC and CQC, we need the following version of the DNS schema, which is extended with a clause handling implication, something that is not needed when one has the $\bot_E$ rule. We denote by $\neg_T A$ the formula $A\Rightarrow T$; when it is clear from the context, we omit the subscript $T$ from $\neg_T$.



\begin{definition} The \emph{Double Negation Shift for $T$} (\dnst) is the following generalisation of the minimal-predicate-logic version of the usual DNS schema, extended with a clause handling implication:
  \begin{align}
    \tag{\dnstforall} \forall x. \neg_T\neg_T A(x) & \Rightarrow \neg_T\neg_T \left(\forall x. A(x)\right)\\
    \tag{\dnstimpl} \left(A\Rightarrow\neg_T\neg_T B\right) & \Rightarrow \neg_T\neg_T\left(A\Rightarrow B\right)
  \end{align}
\end{definition}


% ...


The following proposition is given for IQC as Exercise 2.3.3 of \cite{TroelstraVD1}, we give the proof here to emphasise the role of \dnstimpl\ when $\bot_E$ is not present.
\begin{proposition}\label{dneg} $\text{\dnst}\vdash^m \neg_T\neg_T A \Leftrightarrow A^T$.
\end{proposition}
\begin{proof} Induction on the complexity of $A$. When $A$ is atomic, $A^T=\neg\neg A$.
  \begin{itemize}
  \item[$(\wedge)$] Both directions are via the proof term
    \begin{align*}
      \lambda c.\lambda k. \text{IH}_A & \left(\lambda k'. c \left(\lambda d. k'\left(\pi_1 d\right)\right)\right)\\
        &\left(\lambda a. \text{IH}_B \left(\lambda k'. c\left(\lambda d. k'\left(\pi_2 d\right)\right)\right) \left(\lambda b. k\left(a,b\right)\right)\right).
    \end{align*}
  \item[$(\vee)$] Both directions are via the proof term
    \begin{align*}
      & \lambda a.\lambda k. a \left(\lambda c.\right. \\
      & \left.\caseof{c}{a_1}{\text{IH}_A \left(\lambda l.l a_1\right) \left(\lambda b.k\left(\iota_1 b\right)\right)}{a_2}{\text{IH}_B \left(\lambda l.l a_2\right) \left(\lambda b.k\left(\iota_2 b\right)\right)}\right)
    \end{align*}
  \item[$(\exists)$] Analogous to case $(\vee)$.
  \item[$(\Rightarrow)$] 
    In this case it is crucial to use \dnstimpl, since in minimal logic we do not have $\bot_E$.
        
    The direction left-to-right is via the proof term
    \begin{align*}
      \lambda c.\lambda k. \text{IH}^\leftarrow_A & \left(\lambda k'. \text{\dnstimpl} \left(\lambda a.\lambda k''.k' a\right) k\right)\\
      & \left(\lambda a.\text{IH}^\rightarrow_B \left(\lambda k'. c \left(\lambda f. k' \left(f a\right)\right)\right) \left(\lambda b.k\left(\lambda a'. \lambda k'. k' b\right)\right)\right).
    \end{align*}

    The direction right-to-left is via the proof term
    \begin{align*}
      \lambda c.\lambda k. \text{IH}^\rightarrow_A & \left(\lambda k'. \text{\dnstimpl} \left(\lambda a.\lambda k''. k' a\right) k \right)\\
      & \left(\lambda a.\text{IH}^\leftarrow_B \left(\lambda k'. c \left(\lambda f. f a k'\right)\right) \left(\lambda b.k\left(\lambda a'. b\right)\right)\right).
    \end{align*}

    The arrows in the superscript of ``IH'' determine the direction in which the induction hypotheses are used.
  \item[$(\forall)$] We have:
    \begin{align*}
      (\forall x A(x))^T = &\neg\neg(\forall x A^T(x)) \overset{\text{IH}}{\leftrightarrow} \neg\neg(\forall x \neg\neg A(x)) \\
      \overset{\text{\dnstforall}}{\leftrightarrow} & \neg\neg\neg\neg\forall x A(x) \leftrightarrow \neg\neg\forall x A(x)
    \end{align*}
  \end{itemize}
\end{proof}

\begin{lemma}\label{ctom} $\Gamma \vdash^c A$ if and only if $\Gamma_T\vdash^m A^T$.
\end{lemma}
\begin{proof} The direction right-to-left follows from the previous proposition, because DNS is a classical theorem. The other direction is by induction on the derivation of $\Gamma\vdash^c A$. Actually, we can use the translation table of the proof of Theorem~\ref{equiconsistency} to treat all cases, except for the $\neg\neg_E$ rule which was not covered by the translation. We remark that there is no need to translate the $\bot_E$ rule, since it comes for free in classical logic -- it is derivable from the $\neg\neg_E$ rule. 

To show that $\Gamma_T\vdash^m A^T$ follows from $\Gamma_T\vdash^m (\neg\neg A)^T$, we use the fact that $\vdash^m \neg\neg (T_T) \leftrightarrow T$:
\begin{align*}    
  (\neg\neg A)^T &= ((A\Rightarrow T)\Rightarrow T)^T = \neg\neg((A_T\Rightarrow \neg\neg T)\Rightarrow \neg\neg T) \\
  & \Leftrightarrow \neg\neg ((A_T\Rightarrow T)\Rightarrow T) = \neg\neg\neg\neg A_T \Leftrightarrow \neg\neg A_T = A^T.
  \end{align*}
\end{proof}

We proved the following relationships for the provability of an arbitrary formula $A$ in \mqcplus, MQC, and CQC:
\[
\xymatrix{
  \vdash^+ A \ar[r]^{\ref{equiconsistency}} & \vdash^m A^T \ar@{->}[d]^{\ref{dneg}}  \ar@{<->}[r]^{\ref{ctom}} &  \vdash^c A \\
  \vdash^+\neg\neg A &  \text{\dnst}\vdash^m \neg\neg A \ar[l] }
\]

\begin{corollary}For any formula $A$, we have the following diagram:
\[
\xymatrix{
  \vdash^+ \neg\neg A \ar[r]^{\ref{equiconsistency}} & \vdash^m (\neg\neg A)^T \ar@{->}[d]^{\ref{dneg}} \ar@{<->}[r]^{\ref{ctom}} &  \vdash^c A \\
  \text{\dnst}\vdash^m\neg\neg A \ar[u] &  \text{\dnst}\vdash^m \neg\neg\neg\neg A \ar[l]}
\]
In particular, the statement $\vdash^+\neg\neg A \longleftrightarrow ~ \vdash^c A$ represents an extension of Glivenko's theorem \cite{Glivenko1929,sep-intuitionistic-logic-development,wiki:Glivenko} to predicate logic.
\end{corollary}


\section{Properties}
\label{delcont_srprogress}

In this section we will prove that \mqcplus\ has the Normalisation, Disjunction, and Existence Properties, by proving properties of the reduction relation on proof terms.

\begin{lemma}[Annotation Weakening] \label{lem:weaken}If $\Gamma\vdash p:A$, then $\Gamma\vdash_T p:A$ for any $T$.
\end{lemma}
\begin{proof}A simple induction on the derivation.\end{proof}

\begin{lemma}[Substitutions] \label{lem:subst}The following hold:
  \begin{enumerate}
  \item If $\Gamma, a:A\vdash_\diamond p:B$ and $\Gamma\vdash_\diamond
    q:A$, then $\Gamma\vdash_\diamond p\{q/a\}:B$.
  \item If $\Gamma\vdash_\diamond p:B(x)$, where $x$ is fresh, and $t$ is a
    closed term, then $\Gamma\vdash_\diamond p\{t/x\}:B(t)$.
  \end{enumerate}
\end{lemma}
\begin{proof}The proof is standard, by induction on the derivation (see for example \cite{SorensenU}). The new rules $\mathcal{S}$ and $\#$ pose no problems, since we can use the identities $(\reset{p})\{q/a\} = \reset{(p\{q/a\})}$ and $(\shift{k}{p})\{q/a\} = \shift{k}{(p\{q/a\})}$ when $k$ is fresh.
\end{proof}

\begin{lemma}[Decomposition]\label{lem:decomp}
  If $\Gamma\vdash_T P[\shift{k}{p}] : B$, then there is a
  formula $A$ and derivations $\Gamma, k:A\Rightarrow
  T\vdash_T p:T$ and $\Gamma, a:A\vdash_T P[a] : B$.
\end{lemma}
\begin{proof}
  The proof is by induction on the derivation. We only need to consider the rules that can generate a pure evaluation context of the required form. Of the rules that we consider, for the intuitionistic rules, the proof is simply by using the induction hypothesis, as shown below for the $\wedge_I$ rule; and the only non-intuitionistic rule to consider is $\mathcal{S}$, because $\#$ does not generate a pure evaluation context.
  \begin{itemize}
  \item For $\wedge_I$, there are two cases to consider, depending on whether the pure evaluation context is $(P[\shift{k}{p}],q)$ or $(V,P[\shift{k}{p}])$, but the proofs are analogous. Let the last rule in the derivation be:
    \begin{prooftree}
      \axc{$\Gamma\vdash_T P[\shift{k}{p}] : B_1$}
      \axc{$\Gamma\vdash_T q : B_2$}
      \bic{$\Gamma\vdash_T (P[\shift{k}{p}],q) : B_1\wedge B_2$}{}
    \end{prooftree}
    The induction hypothesis gives us a formula $A_1$ and two derivations, $\Gamma,k:A_1\Rightarrow T\vdash_T p:T$ and $\Gamma,a:A_1\vdash_T P[a]:B_1$, from which the goal follows by choosing $A:=A_1$.
  \item For $\mathcal{S}$, the pure evaluation context must be the empty one, so the last used rule is:
    \begin{prooftree}
      \axc{$\Gamma, k:B\Rightarrow T\vdash_T p : T$}
      \uic{$\Gamma\vdash_T [\shift{k}{p}] : B$}{}
    \end{prooftree}
    If we set $A:=B$, the goal follows from the premise of the rule above and, for $\Gamma, a:A\vdash_T [a] : A$, from the \textsc{Ax} rule.
  \end{itemize}
\end{proof}

\begin{lemma}[Annotation Strengthening]\label{lem:strength}
  $\Gamma\vdash_S V:T \longrightarrow \Gamma\vdash V:T$
\end{lemma}
\begin{proof}The proof is by induction on the derivation and very simple. We only need to consider the intuitionistic rules that introduce a value and that prove a $\Sigma$-formula, that is, the rules \textsc{Ax}, $\wedge_I$, $\vee^1_I$, $\vee^2_I$, and $\exists_I$. $\mathcal{S}$ and $\#$ do not introduce a value.
\end{proof}

\begin{theorem}[Subject Reduction]\label{thm:sr} If
  $\Gamma\vdash_\diamond p:A$ and $p\to q$, then $\Gamma\vdash_\diamond
  q:A$.
\end{theorem}
\begin{proof} The proof is by induction on the derivation and is standard (see for example \cite{SorensenU}), by using Substitutions Lemma \ref{lem:subst} and Decomposition Lemma \ref{lem:decomp}. Below, we consider the new rules and, for illustration, one of the intuitionistic rules.
  \begin{itemize}
  \item[($\#$)] We have $\Gamma\vdash_\diamond \reset{p}$ and $\reset{p}\to q$ for some $q$. We look at three possible cases, because there are three rules for reducing a term of form $\reset{p}$. If $q\equiv \reset{q'}$ and the reduction was by the congruence rule, we have $p\to q'$; now use IH and the $\#$ rule to finish the proof. If $p$ is a value and $q\equiv p$, then $\Gamma\vdash_T q:T$; now use Strengthening Lemma \ref{lem:strength} to conclude $\Gamma\vdash q:T$. The third case is when $p\equiv P[\shift{k}{p'}]$ and $q\equiv \reset{p'\{(\lambda a.\reset{P[a]}) / k\}}$ -- then, the proof is by combining lemmas \ref{lem:subst} and \ref{lem:decomp}.
  \item[($\mathcal{S}$)] This case is impossible, since there are no rules for reducing a term of form $\shift{k}{p}$ on its own, and the set of evaluation contexts does not include a clause for $\shift{k}{[~]}.$
  \item[$(\wedge_E^1)$] We have $\Gamma\vdash_\diamond p:A\wedge B$, $\Gamma\vdash_\diamond \pi_1 p:A$, and $\pi_1 p \to q$. If the reduction was by the congruence rule, then $q\equiv \pi_1 q'$ for some $q'$, and we can use IH. Otherwise, $p\equiv (V_1,V_2)$ and $q\equiv V_1$, and $\Gamma\vdash_\diamond p:A\wedge B$ must have been proved by the $\wedge_I$ rule, which is enough.
  \end{itemize}
\end{proof}

While the last theorem shows that reducing a proof term does not change its logical specification, the next one shows that a proof term which is not in normal form does not get ``stuck''.

\begin{theorem}[Progress]\label{progress} If $\vdash_\diamond p:A$, $p$ is not a value, and $p$ is not of form $P[\shift{k}{p'}]$, then $p$ reduces in one step to some proof term $r$.
\end{theorem}
\begin{proof}
  By induction on the derivation. The cases \textsc{Ax}, $(\Rightarrow_I)$, and $(\forall_I)$ introduce a value, while the case $(\mathcal{S})$ introduces a $\shift{k}{p}$ term, so they are impossible.
  \begin{itemize}
  \item[$(\wedge_I)$] We have that $\vdash_\diamond (p,q):A\wedge B$ and $(p,q)$ is neither a value nor of form $P[\shift{k}{p'}]$. From $(p,q)\neq V$, we have that at least one of $p$ and $q$ is not a value. From $(p,q)\neq P[\shift{k}{p'}]$, we have that $p\neq P'[Sk.r]$, and $p\neq V'$ or $q\neq P'[\shift{k}{r}]$.

    If $p$ is not a value, since it is neither of form $P'[\shift{k}{r}]$, by IH, $p\to r$ for some $r$, hence $(p,q)\to (r,q)$.

    If $p$ is a value $V'$, then $q$ is not a value and it is not of form $P'[\shift{k}{r}]$, and, by IH, $q\to r$, therefore $(V',q)\to (V',r)$.
  \item[$(\wedge_E^1)$] We have that $\vdash_\diamond \pi_1 p : A$ and that $\pi_1 p$, hence $p$ itself, is not of form $P[\shift{k}{p'}]$. If $p$ is a value, then it must be a pair $(V_1,V_2)$, so $\pi_1 (V_1,V_2) \to V_1$. If $p$ is not a value, we can use IH to obtain $\pi_1 p \to \pi_1 r$ for some $r$.
  \item[$(\vee^1_I)$] From $ \vdash_\diamond \iota_1 p:A\vee B$ and $\iota_1 p$ a non-value and not of form $P[\shift{k}{p'}]$, we have that $p$ is not a value and not of that form, so we use IH to obtain an $r$ such that $p \to r$, hence $\iota_1p \to \iota_1r$.
  \item[$(\vee_E)$] We have $\vdash_\diamond \caseof{p}{a_1}{p_1}{a_2}{p_2} : C$. If $p$ is a value, then it is of form $\iota_i V$, therefore $\caseof{\iota_i V}{a_1}{p_1}{a_2}{p_2} \to p_i\left\{V/a_i\right\}$. If $p$ is of form $P[\shift{k}{p'}]$, then so is $\caseof{p}{a_1}{p_1}{a_2}{p_2}$. Otherwise, we use IH to obtain an $r$ such that $\caseof{p}{a_1}{p_1}{a_2}{p_2} \to \caseof{r}{a_1}{p_1}{a_2}{p_2}$.
  \item[$(\Rightarrow_E)$] From $p q\neq P[\shift{k}{r}]$, we have that $p\neq P'[\shift{k}{r}]$, and $p\neq V'$ or $q\neq P'[\shift{k}{r}]$.

    If $p$ is not a value, since it is also not of form $P'[\shift{k}{r}]$, we can use IH.

    If $p$ is a value, it is of form $\lambda a.p'$. If $q$ is a value $V$, then $p q\to p'\{V/a\}$. If $q$ is not a value, it can not be of form $P'[\shift{k}{r}]$, because $p$ is a value; then, we can use IH on $q$.
  \item[$(\forall_E)$] We have $\vdash_\diamond p t:A(t)$. If $p$ is of form $P[\shift{k}{p'}]$, then so is $p t$. If $p$ is a value, then it is of form $\lambda x. r$, hence $(\lambda x. r) t \to r\{t/x\}$. Otherwise, by IH, $p\to r$ for some $r$, so $p t \to r t$.
  \item[$(\exists_I)$] From $ \vdash_\diamond (t,p):A(t)$ and $(t,p)$ a non-value and not of form $P[\shift{k}{p'}]$, we have that $p$ is not a value and not of that form, so we use IH to obtain an $r$ such that $(t,p) \to (t,r)$.
  \item[$(\exists_E)$] We have $\vdash_\diamond \destas{p}{x}{a}{q} : C$. If $p$ is a value, then it is of form $(t,V)$, therefore $\destas{(t,V)}{x}{a}{q} \to q\left\{t/x\right\}\left\{V/a\right\}$.  If $p$ is of form $P[\shift{k}{p'}]$, then so is $\destas{p}{x}{a}{q}$. Otherwise, we use IH to obtain an $r$ such that $\destas{p}{x}{a}{q} \to \destas{r}{x}{a}{q}$.
  \item[$(\#)$] We have $\vdash_\diamond \reset{p}:T$. If $p$ is a value, then $\reset{p} \to p$.  If $p\equiv P[\shift{k}{p'}]$, then $\reset{p}\to \reset{p'\{\lambda a.\reset{P[a]}/k\}}$. If $p$ is neither a value nor of form $P[\shift{k}{p'}]$, by IH, $p\to p'$, so $\reset{p}\to \reset{p'}$.
  \end{itemize}
\end{proof}

\begin{corollary}[Normalisation] For every closed proof term $p_0$, such that $\emptyset \vdash^+ p_0:A$, there is a finite reduction path $p_0\to p_1\to \ldots \to p_n$ ending with a value $p_n$.
\end{corollary}
\begin{proof}
This is a consequence of Subject Reduction and Progress, because a derivation tree $\emptyset \vdash^+ p_0:A$, with no annotations under the turnstile, can not reduce to the form $P[\shift{k}{p}]$. That the reduction path has finite length follows from Theorem~4 of \cite{AsaiK2007,AsaiK2007TR}. 
\end{proof}

\begin{remark}
When proving Normalisation of a variant of $\lambda$-calculus, it is customary to distinguish \emph{weak} (there exists a terminating reduction sequence) from \emph{strong} normalisation (all reduction sequences terminate). There is no need to make that distinction in the present case, because the reduction system is of the type known as \emph{weak head reduction}, which only permits one possible reduction sequence.
\end{remark}

\begin{corollary}[Disjunction and Existence Properties]
  If $\emptyset \vdash^+ A\vee B$, then $\emptyset \vdash^+ A$ or $\emptyset \vdash^+ B$.  If $\emptyset \vdash^+ \exists x A(x)$, then there exists a closed term $t$ such that $\emptyset \vdash^+ A(t)$.
\end{corollary}
\begin{proof}
  Let $\emptyset \vdash^+ p:A\vee B$. By Normalisation and Subject Reduction, for some $V$, $p \to \cdots \to V$ and $\emptyset \vdash^+ V:A\vee B$. Since $V$ is a value, $V$ must be of form $\iota_1 V'$ or $\iota_2 V'$, therefore either $\emptyset \vdash^+ V':A$ or $\emptyset \vdash^+ V':B$. The case for ``$\exists$'' is analogous.
\end{proof}


\section{Related and future work}
\label{delcont_future}
\subsection{Double-negation Shift}
The first use of a schema equivalent to DNS appears to be in modal logic, by Barcan \cite{Barcan1946,Barcan1946a,Fitting1997}, who introduced what is today known as Barcan's formula,
\[
\forall x \Box A(x) \to \Box\forall x A(x),
\]
or, equivalently,
\[
\Diamond\exists x A(x) \to \exists\Diamond A(x).
\]

Veldman kindly pointed to us that DNS is also known as Kuroda's Conjecture \cite{Kuroda1951}. In \cite{Kripke1965}, Kripke showed that Kuroda's Conjecture and Markov's Principle are underivable in intuitionistic logic. (however, see also \cite{Kreisel1970} for criticism of Kripke's argument) 

In \cite[Section 2.11]{Kreisel1957}, Kreisel used the principle
\[\tag{GMP}
\neg\forall n A(n) \Rightarrow \exists n \neg A(n),
\]
for $A(n)$ an arbitrary formula, to deal with implication while giving a translation of formulae of Analysis into functionals of finite type. In 
\cite{Oliva2006}, Oliva calls this principle the Generalised Markov's Principle (GMP) and remarks that $\text{HA}^\omega \vdash \text{DNS} \leftrightarrow \neg\neg \text{GMP}$. Kreisel does not give a justification of GMP in his paper.

The term ``double negation shift'' appears for the first time in \cite{Spector} to denote the formula
\[\tag{DNS}
\forall n \neg\neg A(n) \Rightarrow \neg\neg\forall n A(n).
\]
There, Spector builds upon previous works of G\"odel \cite{Godel1941, Godel1958, Godel1972}, namely he realises DNS by adding the schema of bar recursion to Gödel's system T. The name ``bar recursion'' comes from the Bar Principle of Brouwer which is used in justifying it. However, Spector attaches no particular interest to the DNS schema itself; he writes:
\begin{quote}
  The schema [DNS] is chosen not because we believe it is of intuitionistic significance, but to provide a formal system in which classical analysis is easily interpreted, and whose logical basis is intuitionistic. \cite{Spector}
\end{quote}

We treat DNS at the level of predicate logic, not of arithmetic, an important change in status that we plan to investigate in future.

\subsection{Negative translation of Countable Choice} The Axiom of Countable Choice,
\[\tag{\cchoice}
\forall x^0\exists y^\rho A(x,y) \Rightarrow \exists f^{0\to \rho}\forall x^0 A(x,f(x)),
\]
is a formula schema of $\text{HA}^\omega$, Heyting Arithmetic in all finite types. The type $0$ stands for the set of natural numbers $\mathbb{N}$, the type $1=0\to 0$ stands for the functions $\mathbb{N}\to\mathbb{N}$, $2$ stands for the functionals $(0\to 0)\to 0$, and so on; $\rho$ is a type variable. 

Spector showed that Kuroda's \cite[p.163]{Kohlenbach} negative translation, $\neg\neg(\text{\cchoice}^*)$, of \cchoice,

\[\tag{\cchoicen}
\forall x^0\neg\neg\exists y^\rho A^*(x,y) \Rightarrow \exists f^{0\to \rho}\forall x^0 \neg\neg A^*(x,f(x)),
\]
is provable from DNS and the intuitionistic \cchoice. Since \cchoice\ is realisable in HA$^\omega$, and DNS is realisable by bar recursion, so is \cchoicen. His approach was extended to the Axiom of Dependent Choice (\dchoice) by Luckhardt \cite{luckhardt} and Howard \cite{howard}. In more recent years, Kohlenbach, Berger, and Oliva have given their own versions of bar recursion (see \cite{BergerO2005} for a comparison).


Since we treat DNS at the level of pure logic, without considering arithmetic axioms, we are only able to give an \emph{open} proof term deriving the negative translation of \cchoice,
\[\tag{\cchoice$_T$}
\forall x^0\neg_T\neg_T\exists y^\rho A_T(x,y) \Rightarrow \neg_T\neg_T\exists f^{0\to \rho}\forall x^0 \neg_T\neg_T A_T(x,f(x)).
\]

Given a variable $c$ to denote a proof of the intuitionistic \cchoice, we can use a proof term similar to the one of \dnst\ for deriving the above schema:
\[
\lambda a. \lambda k. \reset{k (c (\lambda x. \shift{k'}{a x (\lambda d. k' (\nu d))}))},
\]
where $\nu$ is a proof term for $\exists y A_T(x,y) \Rightarrow \exists y A^T(x,y)$.

The proof term being open means that we can not immediately use it for computation. We would have to either develop a realisability interpretation for \mqcplus, or add delimited control operators to an intuitionistic system with strong existential quantifiers, like Martin-Löf's type theory, which can derive \cchoice.

\subsection{Herbelin's calculus for Markov Principle}
\label{delcont_future_markov}
In \cite{Herbelin2010}, Herbelin presented IQC$_\text{MP}$, an intuitionistic predicate logic that can derive the pure predicate-logical version of Markov's Principle. Our \mqcplus\ has been developed starting from his calculus. There are two important differences between the two.

First, derivations of IQC$_\text{MP}$ are annotated by a \emph{context} of $\Sigma$-formulae, not just one formula. This permits to have a derivation which uses multiple and different instances of Markov's Principle. Had we had context-annotations as well, it would have been possible to have the following characterisation of provability of $\Sigma$-formulae $S$:
\[
\xymatrix{
  \vdash^+ S \ar[r]^{\ref{equiconsistency}} & \vdash^i S^\bot \ar@{=}[d]^{\text{by def. of $(\cdot)^\bot$}}  \ar@{<->}[r]^{\ref{ctom}} &  \vdash^c S \\
  \text{MP}\vdash^i S \ar[u] &  \vdash^i \neg\neg S \ar[l]}
\]
Proving that such a context-annotated version of \mqcplus\ satisfies the analogues of the properties proven in Section~\ref{delcont_srprogress} remains future work.

Aside from that, the typing and the reduction rules for delimited control operators of IQC$_\text{MP}$ are a restriction of those for \mqcplus.  Consider the typing rules:

\begin{tabular}{ m{4cm} m{5cm} }
  \begin{prooftree}
    \axc{$\Gamma\vdash_{\alpha:T,\Delta} p:T$}
    \uic{$\Gamma\vdash_\Delta \mathsf{catch}_\alpha p:T$}{\textsc{Catch}}
  \end{prooftree}
  &
  \begin{prooftree}
    \axc{$\Gamma \vdash_\Delta p:T \quad (\alpha:T)\in\Delta $}
    \uic{$\Gamma\vdash_\Delta \mathsf{throw}_\alpha p :A$}{\textsc{Throw}}
  \end{prooftree}\\
\end{tabular}

While $\mathsf{catch}$ is just $\#$, the proof term $\mathsf{throw}~p$ is a particular case of $\shift{k}{p}$ that does not use the continuation variable $k$ inside $p$, something already seen with the proof term deriving MP of Example \ref{example_derivations_terms}.

\subsection{Other studies of delimited control operators}

Delimited control operators have been studied in the Theoretical Computer Science literature quite extensively in the past twenty years, since their appearance \cite{Landin1965,Thielecke1998,FelleisenFKD1986}. We mention some of the works that pertain to Logic.

The original typing system for shift/reset of Danvy and Filinski from \cite{DanvyF1989} is a so-called ``type-and-effect'' system: implication is a quaternary not a binary connective, and is as such difficult to understand in traditional logical terms. A proof of Subject Reduction, Progress, and Normalisation of this system appears in \cite{AsaiK2007,AsaiK2007TR}.

A typing system which is a specialisation of Danvy and Filinski's, but again has a ternary implication connective, appears in \cite{Murthy1992}. 

Ariola, Herbelin and Sabry \cite{AriolaHS2007} decompose shift/reset in their own calculus and prove, besides other things, the normalisation property for a typing system where reset is applied at atomic types.

There are a number of works connecting delimited control operators to sub-structural classical logic \cite{KiselyovS07,Zeilberger2010,Munch2011}. Our contribution differs in two respects: we identify delimited control operators as giving rise to an \textit{a priori} constructive logic, rather than classical logic which is only constructive \emph{a posteriori} for certain classes of formulae; and we are connecting delimited control operators with known extra-intuitionistic axioms, rather than analysing the sub-structural properties of the derivation system rules themselves. 

\subsection{The meaning of DNS in the presence of common axioms}

In this paper we were dealing with a purely predicate-logical version of DNS.  Further examination is necessary on how DNS interacts with common logical axioms, as Wim Veldman kindly warned us. For instance, DNS is false in some uncountable models: for example, it contradicts the continuity principle proposed by Brouwer.


\section*{Acknowledgements} I would like to thank my thesis supervisor Hugo Herbelin for commenting on an earlier draft of this paper, as well as for many inspiring discussions.

\bibliographystyle{plain}
\bibliography{article}

\end{document}